\documentclass[12pt,amsymb,fullpage]{amsart}
\usepackage{amssymb,amscd,pstricks}
\newtheorem{Theorem}{Theorem}[section]

\newtheorem{Lemma}[Theorem]{Lemma}

\newtheorem{theorem}{Theorem}[section]
\newtheorem{defn}[theorem]{Definition}

\newtheorem{lemma}[theorem]{Lemma}

\newtheorem{eple}[theorem]{Example}
\newtheorem{rmk}[theorem]{Remarks}
\newtheorem{dsc}[theorem]{Discussion}
\newtheorem{nota}[theorem]{Notation}

\newsavebox{\indbin}
\savebox{\indbin}{\begin{picture}(0,0)
\newlength{\gnu}
\settowidth{\gnu}{$\smile$} \setlength{\unitlength}{.5\gnu}
\put(-1,-.65){$\smile$} \put(-.25,.1){$|$}
\end{picture}}

\newcommand{\be}{\begin{enumerate}}
\newcommand{\bd}{\begin{defn}}
\newcommand{\bt}{\begin{theorem}}
\newcommand{\bl}{\begin{lemma}}
\newcommand{\ee}{\end{enumerate}}
\newcommand{\ed}{\end{defn}}
\newcommand{\et}{\end{theorem}}
\newcommand{\el}{\end{lemma}}

\begin{document}
\title{Infintesimals in a Recursively Enumerable Prime Model}
\author{Tristram de Piro}
\address{Mathematics Department, The University of Camerino, Camerino, Italy}
 \email{tristam.depiro@unicam.it}
\thanks{The author was supported by a Marie Curie research fellowship}
\begin{abstract}
Using methods developed by Robinson, we find a complete theory
suitable for a first order description of infintesimal
neighborhoods. We use this to construct a specialisation having
universal properties and to find a recursively enumerable model in
which the algebraic version of Bezout's theorem is provable by
non-standard methods.
\end{abstract}

\maketitle

\begin{section}{Specialisations and Valuations}

Let $L$ and $K$ be fields with an imbedding $i:L^{*}\rightarrow
K^{*}$. In the case when $L$ and $K$ have the same characteristic,
we will consider $L$ as a subfield of $K$, otherwise we will by
some abuse of notation refer to the embedded set $i(L^{*})\cup
\{0\}$ as L. Let $P(K)=\bigcup_{n\geq 1}P^{n}(K)$ and
$P(L)=\bigcup_{n\geq 1}P^{n}(L)$. By a closed algebraic subvariety
of $P^{n}(K)$, we mean a set $W(K)$ where $W$ is defined by
homogeneous polynomial equations with coefficients in $K$. We say
that $W(K)$ is defined over $L$ if we can take the coefficients to
lie in $L$. Let $W_{n}^{m}(K)$ denote the $m'th$ Cartesian product
of $P^{n}(K)$. By a closed algebraic subvariety of $W_{n}^{m}(K)$,
we mean a set $W(K)$ defined by multi-homogeneous polynomial
equations with coefficients in $K$, similarly we can make sense of
the notion of being defined over $L$. Note that if $K$ is not
algebraically closed, it is not necessarily true that the
projection maps $pr_{k,m}:W_{n}^{k}(K)\rightarrow W_{n}^{m}(K)$
preserve closed
algebraic subvarieties.\\

\begin{defn}

A specialisation is a map
$\pi=\bigcup_{n\geq1}\pi_{n}:P(K)\rightarrow P(L)$, such that
each $\pi_{n}:P^{n}(K)\rightarrow P^{n}(L)$ has the following property;\\

Let $W_{n}^{m}(K)$ denote the $m$'th Cartesian product of
$P^{n}(K)$. Then, if $V\subset W_{n}^{m}(K)$ is a closed algebraic
subvariety defined over $L$ and $\bar a$ is an $m$-tuple of
elements from $W_{n}(K)$, such that $V(\bar a)$ holds, then
$V(\pi_{n}(\bar a))$ holds as well.\\

The following compatibility requirement must also hold between the
$\pi_{n}$;\\

Fix the following chain of embeddings $i_{n}$ of $P^{n}(K)$ and
$P^{n}(L)$ into $P^{n+1}(K)$ and $P^{n+1}(L)$ for $n\geq 1$.\\

$i_{n}:[x_{0}:\ldots:x_{n}]\mapsto [x_{0}:\ldots:x_{n}:0]$.\\

Then we require that $\pi_{n+1}\circ i_{n}=i_{n+1}\circ
\pi_{n}$.\\

\end{defn}

\begin{defn}

A Krull valuation $v$ is a map $v:K\rightarrow \Gamma\cup\infty$
where $\Gamma$ is an ordered abelian group with the following properties;\\

(i). $v(x)=\infty$ iff $x=0$.\\
(ii). $v(xy)=v(x)+v(y)$\\
(iii). $v(x+y)=min\{v(x),v(y)\}$\\

Here, we adopt the convention that $\gamma<\infty$ for
$\gamma\in\Gamma$ and extend $+$ naturally to
$\Gamma\cup\infty$.\\

We let $\mathcal O_{v}=\{x\in K:v(x)\geq 0\}$ be the valuation
ring of $v$ and $\mathcal M_{v}=\{x\in K:v(x)>0\}$ the unique
maximal ideal. We also require;\\

(iv). The inclusion $i:L^{*}\cup {0}\rightarrow {\mathcal
O}_{v}^{*}\cup {0}$ maps $L$ isomorphically onto ${\mathcal
O_{v}}/{\mathcal M_{v}}$, the residue field of $v$.\\

\end{defn}

\begin{defn}

We say that two Krull valuations $v_{1}$ and $v_{2}$ are
equivalent, denoted by $v_{1}\sim v_{2}$ if ${\mathcal
O_{v_{1}}}={\mathcal O_{v_{2}}}$.

\end{defn}

\begin{lemma}

$v_{1}$ and $v_{2}$ are equivalent iff there exists
$\Theta:\Gamma_{1}\rightarrow\Gamma_{2}$ such that $\Theta\circ
v_{1}=v_{2}$.

\end{lemma}

In order to see this, define $\Theta(v_{1}(x))=v_{2}(x)$, this is
well defined as if $v_{1}(x)=v_{1}(x')$, then $v_{1}(x/x')=0$,
hence $x/x'$ and $x'/x$ belong to ${\mathcal O_{v_{1}} }$. If
$v_{1}\sim v_{2}$, then $x/x'$ and $x'/x$ belong to ${\mathcal
O_{v_{2}}}$ as well, which gives that $v_{2}(x)=v_{2}(x')$. One
can easily check that $\Theta$ is an isomorphism of ordered
abelian groups as
required.\\

Our main result in this section is the following;\\

\begin{theorem}

Let $X:=\{\pi:P(K)\rightarrow P(L)\}$ be the set of
specialisations and $Y:=\{v/\sim: v:K\rightarrow \Gamma\}$ be the
set of equivalence classes of Krull valuations. Then there exists
a natural bijection between $X$ and $Y$. Specifically, there
exists maps $\Phi$ and $\Psi$;\\

$\Phi:X\rightarrow Y$\\

$\Psi:Y\rightarrow X$\\

with $\Psi\circ\Phi=Id_{X}$ and $\Phi\circ\Psi=Id_{Y}$\\

\end{theorem}
We first show;\\

\begin{theorem}

There exists $\Psi:Y\rightarrow X$

\end{theorem}

\begin{proof}

Let $[v]$ denote a class of Krull valuations on $K$. We define a
specialisation map
$\pi_{[v]}$ as follows;\\

Let $(x_{0}:x_{1}:\ldots:x_{n})$ denote an element of $P^{n}(K)$
written in homogeneous coordinates. For some $\lambda\in K$, the
elements $\{\lambda x_{0},\ldots,\lambda x_{n}\}$ will lie in
${\mathcal O}_{v}$ and not all of them will lie in ${\mathcal
M}_{v}$.  Let $\pi:{\mathcal O}_{v}\rightarrow L$ denote the
unique ring morphism such that $\pi\circ i=Id_{L}$ where $i$ is
the inclusion map from $L$ into ${\mathcal O}_{v}$. Then
$(\pi(\lambda x_{0}):\pi(\lambda x_{1}):\ldots:\pi(\lambda
x_{n}))$ defines an element of $P^{n}(L)$. As is easily checked,
the mapping is independent of the choice of $\lambda$ and depends
only on ${\mathcal O_{v}}$, hence we obtain
$\pi_{n,[v]}:P^{n}(K)\rightarrow P^{n}(L)$. We need to check that
each $\pi_{n,[v]}$ satisfies the property required of a
specialisation. We will just verify this in the case when $m\leq
2$ for each $n\geq 1$, the other cases are straightforward
generalisations;\\

For $m=1$, let $V\subset P^{n}(K)$ be a closed subvariety defined
over $L$, then $V$ is defined by a system of homogeneous equations
in the variables $\{x_{0},\ldots,x_{n}\}$ with coefficients in
$L$. Taking a tuple $\bar a$ belonging to $V$, we can assume that
the elements $\{a_{0},a_{1},\ldots,a_{n}\}$ belong to ${\mathcal
O}_{v}$. Now, using the fact that the residue map $\pi$ is a ring
homomorphism fixing $L$, the reduced elements
$\{\pi(a_{0}),\pi(a_{1}),\ldots,\pi(a_{n})\}$ also satisfy the
same homogeneous equations as required.\\

For the case when $m=2$, we use the Segre embedding which is defined by;\\

$Segre:P^{n}(K)\times P^{n}(K)\rightarrow P^{n(n+2)}(K)$\\

$((x_{0}:\ldots :x_{n}),(y_{0}:\ldots :y_{n}))\mapsto (x_{0}y_{0}:\ldots:
x_{0}y_{n}:x_{1}y_{0}:\ldots :x_{n}y_{n})$\\

The following diagram is easily checked to commute:\\

\begin{eqnarray*}
\begin{CD}
P^{n}(K)\times P^{n}(K)@>Segre>>P^{n(n+2)}(K)\\
@VV\pi_{n,[v]}\times\pi_{n,[v]} V  @VV\pi_{n(n+2),[v]} V\\
P^{n}(L)\times P^{n}(L)@>Segre>>P^{n(n+2)}(L)\\
\end{CD}
\end{eqnarray*}

.\\
Therefore, it is sufficient to prove that the property holds for
$\pi_{n(n+2),[v]}:P^{n(n+2)}(K)\rightarrow P^{n(n+2)}(L)$ when
$m=1$.This is the case covered above.\\

Finally, we need to check the compatibility requirement for the
$\pi_{n,[v]}$, this is a trivial calculation.\\

Denote the specialisation map we have obtained by $\pi_{[v]}$ and
let $\Psi([v])=\pi_{[v]}$.\\

\end{proof}

We now show;\\

\begin{theorem}

There exists $\Phi:X\rightarrow Y$

\end{theorem}

\begin{proof}

Suppose that we are given a specialisation $\pi$. In particular we
have a map $\pi_{1}:P^{1}(K)\rightarrow P^{1}(L)$ satisfying the
requirements
above. We want to show how to recover a Krull valuation on $K$.\\

Let $\gamma:K\rightarrow P^{1}(K)$ be the map $\gamma:k\mapsto
[k:1]$, so $\pi_{1}\circ\gamma:K\rightarrow P^{1}(L)$. Let
$U\subset P^{1}(L)$ be the open subset defined by $P^{1}\setminus
[1:0]$. Let ${\mathcal O}_{K}=(\pi_{1}\circ\gamma)^{-1}(U)$ and
${\mathcal M}_{K}=(\pi_{1}\circ\gamma)^{-1}([0:1])$. We now claim
the following;\\

\begin{lemma}

${\mathcal O}_{K}$ is a subring of $K$ with $Frac({\mathcal
O}_{K})=K$ and ${\mathcal M}_{K}$ is an ideal of ${\mathcal
O}_{K}$.

\end{lemma}

\begin{proof}

 Suppose that $\{x,y\}\subset {\mathcal O}_{K}$, then both
 $\pi_{1}([x:1])$ and $\pi_{1}([y:1])$ are in $U$. Let $C\subset
 P^{1}(K)\times P^{1}(K)\times P^{1}(K)$ be the closed set defined
 in coordinates $([u:v],[w:x],[y:z])$ by the equation $uwz=yvx$.
 As is easily checked, we have that $C([x:1],[y:1],[xy:1])$. By
 the defining property of $\pi_{1}$,
 $C(\pi_{1}([x:1]),\pi_{1}([y:1]),\pi_{1}([xy:1]))$ also holds.
 Therefore, $C([\lambda:1],[\mu:1],[\alpha,\beta])$ where
 $\lambda,\mu,\alpha,\beta$ are in $L$. By definition of $C$, we
 have $\lambda\mu\beta=\alpha$ which forces $\beta\neq 0$. Hence,
 $\pi_{1}([xy:1])\in U$ and therefore $xy\in {\mathcal O}_{K}$.
 Let $D\subset P^{1}(K)\times P^{1}(K)\times P^{1}(K)$ be defined
 using the same choice of coordinates by the equation
 $uxz+wvz=yvx$. Then we have that $D([x:1],[y:1],[x+y:1])$ and
 therefore $D(\pi_{1}([x:1]),\pi_{1}([y:1]),\pi_{1}([x+y:1]))$.
 Again, we must have $D([\lambda:1],[\mu,1],[\delta,\epsilon])$
 where $\lambda,\mu,\delta,\epsilon$ are in $L$. This forces
 $(\lambda+\mu)\epsilon=\delta$ and therefore $\epsilon\neq 0$, so
 $x+y\in {\mathcal O}_{K}$. Clearly, $1\in {\mathcal O}_{K}$ which
 shows that ${\mathcal O}_{K}$ is a subring of $K$ as required.
In order to see that ${\mathcal M}_{K}$ is an ideal of ${\mathcal
 O}_{K}$, let $x\in {\mathcal O}_{K}$ and $y\in {\mathcal M}_{K}$.
 We have that $C([\lambda:1],[0:1],[\alpha,\beta])$ where
 $\pi_{1}([xy:1])=[\alpha,\beta]$. Then
 $\lambda.0.\beta=1.1.\alpha$ forcing $\alpha=0$ and $\beta=1$, so
 $xy\in {\mathcal M}_{K}$. If $x\in {\mathcal M}_{K}$ and $y\in
 {\mathcal M}_{K}$ we obtain $D([0:1],[0:1],[\delta,\epsilon])$
 where $\pi_{1}([x+y:1]=[\delta,\epsilon]$. Then
 $0.1.\beta+1.0.\epsilon=1.1.\delta$, so $\delta=0$ and
 $\epsilon=1$, hence $x+y\in {\mathcal M}_{K}$ as required.
 Finally, we show that $Frac({\mathcal O}_{K})=K$.
 Suppose $x\notin {\mathcal O}_{K}$, then $\pi_{1}([x:1])=[1:0]$.
 We have that $C([x:1],[1/x:1],[1:1])$, hence
 $C([1:0],[\alpha,\beta],[1:1])$ where
 $\pi_{1}([1/x:1])=[\alpha,\beta]$. This forces
 $1.\alpha.1=0.\beta.1$, hence $\alpha=0$ and $\beta=1$. Therefore
 $1/x\in {\mathcal O}_{K}$ as required.
\end{proof}

We now further claim the following;\\

\begin{lemma}

If $\pi_{1}$ is non-trivial, that is $\pi_{1}$ is not a bijection
between $P^{1}(K)$ and $P^{1}(L)$, then ${\mathcal O}_{K}$ is a
proper subring of $K$

\end{lemma}

\begin{proof}

By the same argument as above we have that
$\pi_{1}\circ\gamma(1/{\mathcal M}_{K})=[1:0]$, hence \emph{if}
${\mathcal O}_{K}=K$, using the previous lemma, we must have that
${\mathcal M}_{K}={0}$. If $\pi_{1}$ is non-trivial, we can find
$x\in K$ and $y\in K$ distinct such that
$\pi_{1}([x:1])=\pi_{1}([y:1])$. By the usual arguments, we then
have that $\pi_{1}([x-y:1])=[0:1]$, so $x-y\in {\mathcal M}_{K}$
contradicting the fact that ${\mathcal M}_{K}=\{0\}$.

\end{proof}

We can now construct a Krull valuation on $K$ by a standard
method. Let $\Gamma=K^{*}/{\mathcal O}_{K}^{*}$ and define
$v:K\rightarrow\Gamma$ by $v(x)=x\ mod\ {\mathcal O}_{K}^{*}$ and
$v(0)=\infty$. Define an ordering on the abelian group $\Gamma$ by
declaring $v(x)\leq v(y)$ iff $y/x\in {\mathcal O}_{K}$. This is
well defined as if $v(x)=v(x')$ and $v(y)=v(y')$, then $y'/y$,
$y/y'$,$x/x'$ and $x'/x$ are all in ${\mathcal O}_{K}$. We have
that $y'/x'=y/x.y'/y.x/x'$ and $y/x=y'/x'.y/y'.x'/x$, therefore
$y'/x'\in {\mathcal O}_{K}$ iff $y/x\in {\mathcal O}_{K}$ as
required. Transitivity of the ordering follows from the fact that
${\mathcal O}_{K}$ is a subring of $K$. $\leq$ is a linear
ordering as if $x\in K^{*}$ and $y\in K^{*}$ then either $x/y$ or
$y/x$ lies in ${\mathcal O}_{K}$. Finally, we clearly have that if
$y/x\in{\mathcal O}_{K}$ then $yz/xz\in {\mathcal O}_{K}$, hence
$v(x)\leq v(y)$ implies $v(x)+v(z)\leq v(y)+v(z)$. This turns
$\Gamma$ into an ordered abelian group. Properties $(i)$ and
$(ii)$ of the axioms for a Krull valuation are trivial to check.
Suppose property $(iii)$ fails, then we can find $x,y$ with
$v(x+y)<v(x)$ and $v(x+y)<v(y)$. Therefore $(x+y)/x\notin
{\mathcal O}_{K}$ and $(x+y)/y\notin {\mathcal O}_{K}$. As $1\in
{\mathcal O}_{K}$, we have that $x/y\notin {\mathcal O}_{K}$ and
$y/x\notin {\mathcal O}_{K}$ which is a contradiction. Finally, we
check property $(iv)$. By definition of $\pi_{1}$, we have that
$L^{*}\subset {\mathcal O}_{K}^{*}$, hence $v|L$ is trivial. If
$k\in {\mathcal O}_{K}^{*}$, we can find $l\in L^{*}$ such that
$\pi_{1}([k:1]=[l:1]$, then $\pi_{1}([k-l:1])=[0:1]$ and $k-l\in
{\mathcal M}_{K}$. It follows that $L$ maps onto ${\mathcal
O}_{K}/{\mathcal M}_{K}$, and ${\mathcal O}_{K}/{\mathcal
M}_{K}\cong L$ as required. Denote the valuation we have obtained
by $v_{\pi}$ and set $\Phi(\pi)=[v_{\pi}]$. This ends the proof of Theorem 1.7 \\

\end{proof}

We now complete the proof of Theorem 1.5;\\
\begin{proof}
\emph{$\Phi\circ\Psi=Id_{Y}$}.\\

 Let $[v]$ be a class of Krull
valuations on $K$ with corresponding specialisation $\pi_{[v]}$
provided by $\Psi$. Let $\pi_{1,[v]}$ be the restriction to
$P^{1}(K)$. By definition, if $k\in {\mathcal O}_{v}$ then
$\pi_{1,[v]}([k:1])=[\pi(k),1]$ where $\pi$ is the residue map for
$v$. If $k\notin {\mathcal O}_{v}$, then
$\pi_{1,[v]}([k:1])=[0,1]$, so we see that ${\mathcal O}_{K}$ as
defined above is exactly ${\mathcal O}_{v}$. The valuation
$v_{\pi_{[v]}}$ constructed from $\pi_{[v]}$ therefore has the
same valuation ring ${\mathcal O}_{v}$, so $v\sim v_{\pi_{[v]}}$
which gives the result.\\

\emph{$\Psi\circ\Phi=Id_{X}$}.\\

 Let $\pi$ be a given
specialisation and $[v_{\pi}]$ the corresponding class of Krull
valuations. Let $\pi_{1}$ be the restriction of $\pi$ to
$P^{1}(K)$ and $\pi_{1,v_{\pi}}$ the specialisation constructed
from $v_{\pi}$ restricted to $P^{1}(K)$. We have;\\

(i). $\pi_{1,v_{\pi}}([k:1])=[0:1]$ iff $v_{\pi}(k)>0$ iff $k\in
{\mathcal M}_{v_{\pi}}$ iff $k\in {\mathcal M}_{K}$ as defined
above iff $\pi_{1}([k:1])=[0:1]$\\

(ii). $\pi_{1,v_{\pi}}([k:1])=[1:0]$ iff $v_{\pi}(k)<0$ iff
$k\notin {\mathcal O}_{v_{\pi}}$  iff $k\notin {\mathcal O}_{K}$
as defined above iff $\pi_{1}([k:1])\notin U$ iff $\pi_{1}([k:1])=[1:0]$\\

(iii). $\pi_{1,v_{\pi}}([1:0])=\pi_{1}([1:0])=[1:0]$ trivially.\\

If $k\in {\mathcal O}_{v_{\pi}}$, then
$\pi_{1,v_{\pi}}([k:1])=[\alpha(k):1]$ where $\alpha$ is the
residue mapping associated to $v_{\pi}$. We also have that
$\pi_{1}([k:1])\in U$, hence as $\pi_{1}$ is a specialisation that
$\pi_{1}([k:1])=[\beta(k):1]$ where $\beta$ is a homomorphism from
${\mathcal O}_{v_{\pi}}$ to $L$. We thus obtain two homomorphisms
$\alpha,\beta:{\mathcal O}_{v_{\pi}}\rightarrow L$ such that (by
(i)) $Ker(\alpha)=Ker(\beta)={\mathcal M}_{v_{\pi}}$ and with the
property that $\alpha\circ i=\beta\circ i=Id_{L}$ where $i$ is the
natural inclusion of $L$ in ${\mathcal O}_{v_{\pi}}$. We thus
obtain the splitting ${\mathcal O}_{v_{\pi}}=L\oplus
Ker(\alpha)=L\oplus Ker(\beta)=L\oplus M$ with
$Ker(\alpha)=Ker(\beta)=M$. Now, using this fact, we can write any
element of ${\mathcal O}_{v_{\pi}}$ uniquely in terms of $L$ and
$M$, hence the corresponding projections $\alpha$ and
$\beta$ are the same.\\

We have shown that $\pi_{1}=\pi_{1,v_{\pi}}$, it remains to check
that $\pi_{n}=\pi_{n,v_{\pi}}$ for all $n\geq 1$. We prove this by
induction on $n$, the case $n=1$ having been established.\\

By the induction hypothesis and the compatibility requirement
between the $\pi_{n}$, for
$\{k_{0},k_{1},\ldots,k_{n}\}\subset {\mathcal O}_{v_{\pi}}$;\\

$\pi_{n+1}([k_{0}:k_{1}:\ldots:k_{n}:0])=[\pi(k_{0}):\pi(k_{1}):
\ldots:\pi(k_{n}):0]$ (*)\\

where $\pi$ is the residue map on ${\mathcal O}_{v_{\pi}}$.\\

Let $C\subset P^{n+1}(K)$ be the closed subvariety defined using
coordinates $[x_{0}:x_{1}:\ldots:x_{n+1}]$ by the equations
$x_{0}=x_{1}=\ldots=x_{n-1}=0$. Then by arguments as above and the
fact that $C$ is preserved by $\pi_{n+1}$, we can find
a Krull valuation $v'$ on $K$ with corresponding residue mapping $\pi'$ such that;\\

$\pi_{n+1}([0:\ldots:0:1:k_{n+1}])=[0:\ldots:0:1:\pi'(k_{n+1})]$\
if $v'(k_{n+1})\geq 0$\\

$\ \ \ \ \ \ \ \ \ \ \ \ \ \ \ \ \ \ \ \ \ \ \ \ \ \ \ \ \ \ \ \
 \ \ \ =[0:\ldots:0:0:1]$\ otherwise (**)\\

Now let $D$ be the closed subvariety of $P^{n+1}(K)$ defined by
the equations $x_{1}=\ldots=x_{n}$ and $x_{0}=x_{n+1}$. Again, we
have that $\pi_{n+1}$ preserves $D$, hence there exists a Krull
valuation $v''$ on $K$ with corresponding residue mapping $\pi''$ such that;\\

$\pi_{n+1}([k:1:\ldots:1:k])=[\pi''(k):1\ldots:1:\pi''(k)]$\ if
$v''(k)\geq 0$\\

$\ \ \ \ \ \ \ \ \ \ \ \ \ \ \ \ \ \ \ \ \ \ \ \ \ \ \ \ \ \ \
=[1:0:\ldots:0:1]$\ otherwise (***)\\

Let $Sum$ be the closed subvariety of $P^{n+1}(K)\times
P^{n+1}(K)\times P^{n+1}(K)$ defined using coordinates
$[x_{0}:x_{1}:\ldots:x_{n+1}]$, $[y_{0}:y_{1}:\ldots:y_{n+1}]$ and
$[z_{0}:z_{1}:\ldots:z_{n+1}]$ by the equations
$x_{0}y_{1}z_{1}+y_{0}x_{n}z_{1}=z_{0}x_{n}y_{1}$ and
$x_{n+1}y_{1}z_{1}+y_{n+1}x_{n}z_{1}=z_{n+1}x_{n}y_{1}$. Then, for
$k\in K$, we have that
$Sum([0:0:\ldots:1:k],[k:1:\ldots:0:0],[k:1:\ldots:1:k])$, hence
by the properties of a specialisation that
$Sum(\pi_{n+1}([0:0:\ldots:1:k]),\pi_{n+1}([k:1:\ldots:0:0]),
\pi_{n+1}([k:1:\ldots:1:k]))$.\\

In the generic case when $v_{\pi}(k),v'(k),v''(k)$ are all
non-negative, we obtain
$Sum([0:0:\ldots:1:\pi'(k)],[\pi(k):1:\ldots:0:0],[\pi''(k):1:\ldots:1:\pi''(k)])$
which gives the relations $0.1.1+\pi(k).1.1=\pi''(k).1.1$ and
$\pi'(k).1.1+0.1.1=\pi''(k).1.1$, so $\pi(k)=\pi'(k)=\pi''(k)$.\\

A simple calculation shows that $v_{\pi}(k)<0$ iff $v'(k)<0$ iff
$v''(k)<0$, hence ${\mathcal O}_{v_{\pi}}={\mathcal
O}_{v'}={\mathcal O}_{v''}$. We have now shown the following
further compatibility between $\pi_{1}$ and $\pi_{n+1}$. Namely;\\

If $\gamma:P^{1}(K)\rightarrow P^{n+1}(K)$ is given by
$\gamma:[x_{0},x_{1}]\mapsto [0:0:\ldots:x_{0}:x_{1}]$ then
$\pi_{n+1}\circ\gamma=\gamma\circ\pi_{1}$.   $(\dag)$\\

Finally, let $Sum'$ be the closed subvariety of $P^{n+1}(K)\times
P^{n+1}(K)\times P^{n+1}(K)$ defined in coordinates
$[x_{0}:\ldots:x_{n+1}],[y_{0}:\ldots:y_{n+1}],[z_{0}:\ldots:z_{n+1}]$
by the $(n+1)$ equations $x_{j}y_{1}z_{1}+
y_{j}x_{n}z_{1}+z_{j}x_{n}y_{1}$ for $j\neq n$. Let
$[k_{0}:\ldots:k_{n+1}]$ be an arbitrary element of $P^{n+1}(K)$.
Without loss of generality, we may assume that
$\{k_{0}:\ldots:k_{n+1}\}\subset {\mathcal O}_{v_{\pi}}$ and that
$k_{n}\in {\mathcal O}_{v_{\pi}}^{*}$. Hence, dividing by $k_{n}$,
the element is of the form $[k_{0}:\ldots:k_{n-1}:1:k_{n+1}]$ with
$\{k_{0},\ldots,k_{n-1},k_{n+1}\}\subset {\mathcal O}_{v_{\pi}}$.
We have that
$Sum'([0:\ldots:0:1:k_{n+1}],[k_{0}:\ldots:k_{n-1}:1:0],[k_{0}:\ldots:k_{n-1}:1:k_{n+1}])$,
hence by specialisation and $(\dag)$,
$Sum'([0:\ldots:0:1:\pi(k_{n+1})],[\pi(k_{0}):\ldots:\pi(k_{n-1}):1:0],
[l_{0}:\ldots:l_{n}:l_{n+1}])$. where
$\{l_{0},\ldots,l_{n+1}\}\subset L$. As is easily checked, the
case when $l_{n}=0$ leads to a contradiction, hence we can assume
that $l_{n}=1$ (multiplying by $1/l_{n}$). Now the equations give
that $l_{j}=\pi(k_{j})$ for $j\neq n$. We have therefore shown
that $\pi_{n+1}=\pi_{n+1,v_{\pi}}$ as required.\\

Theorem 1.5 is now proved.\\
\end{proof}
\end{section}

\begin{section}{A Model Theoretic Language of Specialisations}

We now introduce a model theoretic language which will enable us
to describe specialisations in the context of algebraic geometry.
In this section, we will assume that $K$ and its residue field
have the same characteristic. We will use a many sorted structure
$\{\bigcup S_{n}:n\in {\mathcal N}\}$. Each sort will be the
domain of $P^{n}(K)$ for an algebraically closed field $K$. We fix
an algebraically closed constant field $L$ which we assume to be
countable and let $K$ be some non-trivial extension of $L$, having
the same characteristic. In order to describe algebraic geometry,
we introduce sets of predicates $\{V_{n}^{m}\}$ on the Cartesian
powers $S_{n}^{m}$ to describe closed algebraic subvarieties of
$P^{n}(K)$ defined over $L$. In particular, we have constants to
denote the individual elements of $P^{n}(L)$ on each sort $S_{n}$.
We introduce function symbols $i_{n}:S_{n}\rightarrow S_{n+1}$ to
describe the imbeddings $P^{n}(K)\rightarrow P^{n+1}(K)$ defined
above. Finally, we will have symbols $\{\pi_{n}:n\in {\mathcal
N}\}$ to describe the specialisation map $\pi=\cup_{n\geq
1}\pi_{n}$. Strictly speaking, as $P^{n}(L)$ is not definable,
each $\pi_{n}$ will be a union over $l\in P^{n}(L)$ of unary
predicates defined as $\{x\in P^{n}(K):\pi_{n}(x)=l\}$. We denote
the language $<\{V_{n}^{m}\},i_{n},\pi_{n}>$ by ${\mathcal
L}_{spec}$ and the theory of the structure $<P(K),P(L),\pi>$ in
this language by $T_{spec}$. We denote the theory of the structure
$<P(K),P(L)>$ in the language ${{\mathcal L}_{spec}\setminus
\{\pi_{n}\}}$ by $T_{alg}$. Note that the structure $<K,0,1,+,.>$
is interpretable in the structure $<P(K),P(L)>$ in the language
${{\mathcal L}_{spec}\setminus \{\pi_{n}\}}$ $(*)$. This follows
by noting that the points $[1:0],[0:1]$ and $[1:1]$ are named as
elements in the sort $S_{1}$ and the operations of $+,.$ define
algebraic subvarieties in the sorts $S_{1}^{3}$. The structure
$<L,0,1,+,.>$ is not interpretable but any model of $T_{alg}$ will
contain an isomorphic copy of $P(L)$ as a substructure. It follows
that the models of $T_{alg}$ are exactly of the form $<P(K),P(L)>$
for some algebraically closed field $K$ properly extending $L$
(use the fact that the axiomatisation of $Th(<K,0,1,+,.>)$ can be
interpreted in $T_{alg}$ and the field structure can be related to
the predicates $\{V_{n}^{m}\}$ using the imbeddings $i_{n}$). We
now claim the following;\\

\begin{theorem}

The theory $T_{spec}$ is axiomatised by $T_{axioms}=T_{alg}\cup
\Sigma$ where $\Sigma$ is the set of sentences given by;\\

(i). The mappings $\{\pi_{n}\}$ preserve the predicates
$\{V_{n}^{m}\}$.\\

(ii). The compatibility requirement $\pi_{n+1}\circ
i_{n}=i_{n+1}\circ\pi_{n}$ holds.\\

(see definition 1.1). In particular, $T_{axioms}$ is complete.
Moreover, $T_{axioms}$ is model complete.

\end{theorem}

The proof of this theorem will be based on Theorem 1.5 and the
following result by Robinson, given in \cite{Rob};\\

\begin{theorem}

Let $K$ be an algebraically closed field with a non trivial Krull
valuation $v$ and residue field $l$. Then $T_{K}$ is model
complete in the language ${\mathcal L}_{val}$ and admits
quantifier elimination in the language ${\mathcal L}_{rob}$.
Moreover, the completions of $K$ are determined by the pair
$(char(l),char(K))$, that is $T_{K}\cup\Sigma$ is complete where
$\Sigma$ is the possibly infinite set of sentences specifying the
characteristic of $K$ and $l$.\\

\end{theorem}

Here, by the language ${\mathcal L}_{rob}$ we mean the language of
algebraically closed fields together with a binary predicate
$Div(x,y)$ denoting $v(x)\leq v(y)$. By the language ${\mathcal
L}_{val}$, we mean a $2$-sorted language for the value group and
the field, with the usual language for the field sort and the
language of ordered groups on the group sort. $T_{K}$ is the
theory which asserts that $K$ is an algebraically closed field,
the value group $\Gamma$ is linearly ordered and abelian, the
valuation is non-trivial. For our purposes, we will require a
slightly refined version of this result. Namely, we will fix a set
of constants for an algebraically closed field $L$ which we can
assume to be countable, add to $T_{K}$ the atomic diagram of $L$,
relativized to the field sort, the requirement that $v|L$ is
trivial and $\pi$, the residue mapping, maps $L$ injectively and
homomorphically into the residue field. (Note, the condition that
$L$ maps onto the residue field is not definable and that the
homomorphism requirement ensures that the residue field $l$ and
$K$ have equal characteristic, hence the characteristic of $K$ is
already determined by the characteristic of $L$.) We will denote
the corresponding theory by $T_{K,L}$ and the expanded languages
by ${\mathcal L}_{rob}$ and ${\mathcal L}_{val}$ again. It is no
more difficult to prove that $T_{K,L}$ is model complete,
Robinson's original proof in \cite{Rob} requires the solution of
certain valuation equations in the model $K$ given that these
equations have a solutions in an extension $K'$, it makes no
difference if some of the elements from $K$ are named. In order to
show that $T_{K,L}$ is complete, it is sufficient to
exhibit a prime model of the theory;\\

Case 1. $Char(K,L)=(p,p)$, with $p\neq 0$. Take
$L(\epsilon)^{alg}$ where $\epsilon$ is transcendental over $L$,
define the valuation on $L$ to be zero and extend it to
$L(\epsilon)$ non-trivially using say $v_{ord,\epsilon}$, the
order valuation in $\epsilon$. Take any extension to $L(\epsilon)^{alg}$.\\

Case 2. $Char(K,L)=(0,0)$, define a similar valuation on
$L(\epsilon)^{alg}$.\\

We now show the following lemma;\\

\begin{Lemma}{Amalgamation of Specialisations}\\

Let $(P(K_{1}),P(L),\pi_{1})$ and $(P(K_{2}),P(L),\pi_{2})$ be
models of $T_{axioms}$, then there exists a further model
$(P(K_{3}),P(L),\pi_{3})$ such that;\\

$(P(K_{1}),P(L),\pi_{1})\leq (P(K_{3}),P(L),\pi_{3})$\\

 and\\

$(P(K_{2}),P(L),\pi_{2})\leq (P(K_{3}),P(L),\pi_{3})$\\

\end{Lemma}
\begin{proof}
By Theorem 1.5, we can find Krull valuations $v_{1}$ and $v_{2}$
on $K_{1}$ and $K_{2}$ such that $\pi_{1}=\pi_{v_{1}}$ and
$\pi_{2}=\pi_{v_{2}}$. Using the refined version of Robinson's
completeness result, we can jointly embed $(K_{1},v_{1})$ and
$(K_{2},v_{2})$ over $L$ into $(K_{3},v_{3})$ (*). Let $L'$ be the
residue field of $v_{3}$, then as $K_{3}$ is algebraically closed,
so is $L'$ and extends the residue field $L$ of $v_{1}$ and
$v_{2}$. By standard results, we can construct a Krull valuation
$v$ on $L'$ with residue field $L$, for example use the
construction given in \cite{deP1}. Using Theorem 1.5 again, we can
construct specialisations $(P(K_{3}),P(L'),\pi_{v_{3}})$ and
$(P(L'),P(L),\pi_{v})$, the composition gives a specialisation
$(P(K_{3}),P(L),\pi_{3})$. It remains to see that in fact
$\pi_{3}$ extends the specialisations $\pi_{1}$ and $\pi_{2}$.
This follows from the fact that if $k\in K_{1}$ or $k\in K_{2}$
and there exists $l\in L$ such $v_{1}(k-l)>0$ or $v_{2}(k-l)>0$
then this relation is preserved in the embedding $(*)$. Hence the
specialisation $\pi_{v_{3}}$ already extends the specialisations
$\pi_{1}$ and $\pi_{2}$ of $P(K_{1})$ and $P(K_{2})$ into $P(L)$.
As the specialisation $\pi_{v}$ fixes $L$, this proves the lemma.\\
\end{proof}
\begin{lemma}{Transfer of Formulas}\\

Let $(P(K),P(L),\pi)$ be a specialisation with corresponding
$(K,v)$, then there exists a map;\\

$\sigma:P(K)\rightarrow K^{eq}$\\

$\sigma:{\mathcal L}_{spec}$-formulae $\rightarrow {\mathcal
L}_{val}$-formulae\\

such that for any $\phi(x_{1},\ldots,x_{n})$ which is a ${\mathcal
L}_{spec}$-formula and $(k_{1},\ldots,k_{n})\subset P(K)$;\\

$(P(K),P(L),\pi)\models \phi(k_{1},\ldots,k_{n})$ iff
$(K,v)\models \sigma(\phi)(\sigma(k_{1}),\ldots,\sigma(k_{n}))$\\
\\
.\ \ \ \ \ \ \ \ \ \ \ \ \ \ \ \ \ \ \ \ \ \ \ \ \ \ \ \ \ \ \ \ \
\ \ \ \ \ \ \ \ \ \ \ \ \ \ \ \ \ \ \ \ \ \ \ \ \ \ \ \ \ \ \ \ \
\ \ \ \ \ \ \ \ \ \ \ \ \
$(\dag)$\\
Moreover, the definition of the map is uniform in $K$.\\

\end{lemma}
\begin{proof}
The map $\sigma$ is defined on the sorts $P^{n}(K)$ by sending
$[k_{0},\ldots,k_{n}]$ to $(k_{0},\ldots,k_{n})/\sim_{n}$ where
$\sim_{n}$ is the equivalence relation defined on $K^{n+1}$ from
multiplication by $K^{*}$. Similarly, $\sigma$ maps a variable
from the sort $S_{n}$ to the corresponding variable from the sort
in $K^{eq}$ defined by $\sim_{n}$. A closed algebraic subvariety
in $\{V_{n}^{m}\}$ is defined by a multi-homogeneous equation in
the variables
$\{(x_{01},\ldots,x_{n1}),\ldots,(x_{0m},\ldots,x_{nm})\}$. Let
$C_{n}^{m}$ be the algebraic variety in $K^{m(n+1)}$ defined by
this equation. Then the corresponding formula in $K^{eq}$ is given
by;\\

 $(y_{1},\ldots,y_{m})\in (\sim_{n})^{m}[\exists x_{1}\ldots
x_{m}(C_{n}^{m}(x_{1},\ldots,x_{m})\wedge
\bigwedge_{i=1}^{m}x_{i}/\sim_{n}=y_{i})]$\\

For the inclusion maps $i_{n}$, let us identify each $i_{n}$ with
its graph, then clearly we can define $\sigma$ to map the formula
$i(x)=y$ to a corresponding formula relating the sorts $\sim_{n}$
and $\sim_{n+1}$ in $K^{eq}$.\\

Note that if $l\in P^{n}(L)$ is a constant, then
$\sigma(l)=(l_{0},\ldots,l_{n})/\sim_{n}$ where each $l_{i}$ is a
constant from the atomic diagram of $L$.\\

Finally, let $\pi_{n}:P^{n}(K)\rightarrow P^{n}(L)$ be a
specialisation. Again, let us assume that we can identify
$\pi_{n}$ with its graph. We then have that;\\

$\pi_{n}([x_{0}:\ldots:x_{n}])=[l_{0}:\ldots:l_{n}]$\\

 iff\\

$\exists z\exists z_{0}\ldots\exists
z_{n}((\bigwedge_{i=0}^{n}x_{i}z=l_{i}+z_{i})\wedge(\bigwedge_{i=0}^{n}v(z_{i})>0))$.\\

It is now clear how to define $\sigma(\pi_{n})$ as a union of
formulas in the sort defined by $\sim_{n}$.\\

This completes the definition of $\sigma$, it is clear that the
definition is uniform in $K$ and a straightforward induction on
the length of a formula from ${\mathcal L}_{spec}$ shows that it
has the required property $(\dag)$.\\
\end{proof}

Theorem 2.1 is now a fairly straightforward consequence of the
above lemmas. We first show model completeness. Suppose that we
have models of $T_{axioms}$;\\

$(P(K_{1}),P(L),\pi_{1})\leq (P(K_{2}),P(L),\pi_{2})$\\

By theorem 1.5, we can find Krull valuations $v_{1}$ and $v_{2}$
such that $(K_{1},v_{1})\leq (K_{2},v_{2})$ and
$(K_{1},v_{1}),(K_{2},v_{2})\models T_{K,L}$. By the refined model
completeness result after Theorem 2.2, we have $(K_{1},v_{1})\prec
(K_{2},v_{2})$, hence using Lemma 2.3, we must have that;\\

$(P(K_{1}),P(L),\pi_{1})\prec (P(K_{2}),P(L),\pi_{2})$\\

as required. Completeness now follows directly from Lemma 2.1 and
model completeness. Alternatively, one can exhibit a prime model
of the theory, this is clearly possible by taking the
specialisations corresponding to the prime models of $T_{K,L}$
above.\\

\end{section}

\begin{section}{A First Order Definition of Intersection Multiplicity and
Bezout's Theorem}

We now formulate a non-standard definition of intersection
multiplicity in the language ${\mathcal L}_{spec}$. We will do
this only in the case of projective curves inside $P^{2}(L)$, the
reader may wish to try formulating a corresponding definition in
higher dimensions.\\

Let $C_{1}$ and $C_{2}$ be projective curves of degree d and
degree e in $P^{2}(K)$ defined over $L$. The parameter spaces for
such curves are affine spaces of dimension $(d+1)(d+2)/2$ and
$(e+1)(e+2)/2$ respectively. We can give them a projective
realisation by noting that if $(\bar l)$ is a non-zero vector
defining a curve of degree $d$, then multiplying it by a constant
$\mu$ defines the same curve. Let $P^{d(d+3)/2}(K)$ and
$P^{e(e+3)/2}(K)$ define these spaces which we will denote by
$P_{d}$ and $P_{e}$ for ease of notation. Let $Curve_{d}$ and
$Curve_{e}$ be the closed projective subvarieties of $P_{d}\times
P^{2}(K)$ and $P_{e}\times P^{2}(K)$, defined over the prime
subfield of $L$, such that, for $l\in P_{d}$, the fibre
$Curve_{d}(l)$ defines the corresponding projective curve of
degree $d$ in $P^{2}(K)$. For $l$ in $P^{n}(L)$, we denote its
infintesimal neighborhood ${\mathcal V}_{l}$ to be the inverse
image under the specialisation $\pi_{n}$.\\

 Now suppose that $C_{1}$ and $C_{2}$
(which may not be reduced or irreducible), of degrees $d$ and $e$
respectively, are defined by parameters $l_{1}$ and $l_{2}$ and
intersect at an isolated point $l$ in $P^{2}(L)$. Then we define;\\

$Mult(C_{1},C_{2},l)\geq n$\\

 iff\\

  $\exists x_{1},x_{2}\in {\mathcal
V}_{l_{1}},{\mathcal V}_{l_{2}},\exists_{y_{1}\neq\ldots\neq
y_{n}}\in{\mathcal V}_{l}(\{y_{1},\ldots,y_{n}\}\subset
Curve_{d}(x_{1})\cap
Curve_{e}(x_{2}))$\\

Then define $Mult(C_{1},C_{2},l)=n$\\

 iff\\

  $Mult(C_{1},C_{2},l)\geq
n$ and $\neg Mult(C_{1},C_{2},l)\geq n+1$.\\

Clearly, the statement that $Mult(C_{1},C_{2},l)=n$ naturally
defines a sentence in the language ${\mathcal L}_{spec}$. One
consequence of the completeness result given above is that the
statement "The curves $C_{1}$ and $C_{2}$ intersect with
multiplicity $n$ at $l$" depends only on the theory $T_{axioms}$
and is independent of the particular structure $(P(K),P(L),\pi)$.
In the paper \cite{deP}, we showed that this non-standard
definition of multiplicity is equivalent to the algebraic
definition of multiplicity when computed in the structure
$(P(K_{univ}),P(L),\pi_{univ})$ (see the next section). It
therefore follows that the non-standard definition of multiplicity
is equivalent to the algebraic definition even when computed in a
prime model of $T_{axioms}$ which I will denote by
$(P(K_{prime}),P(L),\pi_{prime})$.\\

We now turn to the statement of Bezout's theorem. In algebraic
language, this says that if projective algebraic curves $C_{1}$
and $C_{2}$ of degree $d$ and degree $e$ in $P^{2}(L)$ intersect
at finitely many points $\{l_{1},\ldots,l_{n}\}$, then;\\

$\sum_{i=1}^{n}I(C_{1},C_{2},l_{i})=de$\\

where $I(C_{1},C_{2},l_{i})$ is the algebraic intersection
multiplicity. The non-standard version of this result can be
formulated in the language ${\mathcal L}_{spec}$ by the
sentence;\\

$Bezout(C_{1},C_{2})\equiv\exists_{m_{1},\ldots,m_{n};m_{1}+\ldots+m_{n}=de}(\bigwedge_{i=1}^{n}Mult(C_{1},C_{2},l_{i})=m_{i})$\\

Again, in the paper \cite{deP}, we proved the \emph{algebraic}
version of Bezout's theorem by non-standard methods in the
structure $(P(K_{univ}),P(L),\pi_{univ})$. It follows that the
sentences $Bezout(C_{1},C_{2})$ are all proved by the theory
$T_{axioms}$ and therefore hold in the structure
$(P(K_{prime}),P(L),\pi_{prime})$ as well. This demonstrates the
fact that we can prove an algebraic statement of Bezout's theorem
using only infintesimals from a straightforward extension of $L$,
namely $L(\epsilon)^{alg}$, in particular in a structure such that
the infintesimal neighborhoods ${\mathcal V}_{l}$ are all
recursively enumerable. This seems to provide some answer to a
general objection concerning the use of infintesimals, originating
in \cite{Berk}. It may also provide an effective alternative
method to compute intersection multiplicities generally in algebraic geometry.\\

\end{section}

\begin{section}{Constructing a Universal Specialisation}

In the papers \cite{deP1} and \cite{deP}, we used the existence of
a specialisation $(P(K_{univ}),P(L),\pi_{univ})$ having the
following "universal" property;\\

 If $L\subset L_{m}$ is an algebraically closed extension of $L$
 with transcendence degree $m$, and $(P(L_{m}),P(L),\pi_{m})$ is a
 specialisation, then there exists an $L$-embedding
 $\alpha_{L}:L_{m}\rightarrow K_{univ}$ with the property that
 $\pi_{univ}\circ\alpha_{L}=\pi_{m}$. (*)\\

 Unfortunately, the construction of $K_{univ}$ was
 flawed. We correct this difficulty here;\\

 Model theoretically, using theorem 2.1, it is easy to show the
 existence of such a structure. Namely, let $(P(K_{univ}),P(L),\pi_{univ})$ be a
 $2^{\omega}$ saturated model of the theory $T_{axioms}$.
 Then, if $L\subset L_{m}$ is an algebraically closed extension of
 $L$ of transcendence degree $m$, clearly $\bigcup_{n\geq
 1}Card(S^{n}(Th({\mathcal M})))\leq 2^{\omega}$, where ${\mathcal
 M}=(P(L_{m}),P(L),\pi_{m})$. This follows as
 $L$ was assumed to be countable. Hence, by elementary model
 theory, there exists an $L$-embedding $\alpha_{L}$ with the
 required properties. For the non-model theorist, we give a more
 algebraic construction, replacing the use of types by an explicit
 amalgamation of the possible valuations;\\

 \begin{proof}
 Suppose, inductively, we have already constructed a specialisation
 $(P(K_{n}),P(L),\pi_{n})$ which has the property $(*)$ for all
 extensions $L\subset L_{m}$ with $L_{m}$ algebraically closed of
 transcendence degree $m\leq n$. We will construct $K_{n+1}$
 having this property for $m\leq n+1$. By theorem 1.5, we can find
 a Krull valuation $v_{n}$ on $K_{n}$ corresponding to the
 specialisation $\pi_{n}$. Let $t$ be a new transcendental
 element. The extensions of $v_{n}$ to $K_{n}(t)$ are completely
 classifiable. In fact, we have the following result in
 \cite{FVK} (Theorem 3.9), we refer the reader to the paper for
 the definition of each family of valuations; \\

The extensions of $v_{n}$ are of the form;\\

(i). $v_{n,a,\gamma}$ where $a\in K_{n}$ and $\gamma$ is an
element of some ordered group extension of $v(K)$.\\

(ii). $v_{n,A}$ where $A$ is a pseudo Cauchy sequence in
$(K_{n},v_{n})$ of transcendental type.\\

Let $I$ be a fixed enumeration of these valuations. Inductively,
we assume that $Card(K_{n})\leq 2^{\omega}$ in which case the
dimension of $v(K_{n})$ as a vector space over ${\mathcal Q}$ has
dimension at most $2^{\omega}$ as well. Clearly then the number of
non-isomorphic (over $K_{n}$) valuations from $(ii)$ is at most
$2^{\omega}$ and the same holds for the valuations obtained from
$(i)$ by noting that the number of order types of $\gamma$ is at
most $2^{\omega}$ (it is easily checked that $2$ new elements of
the value group, $\gamma_{1}$ and $\gamma_{2}$, having the same
order type, define isomorphic valuations in the case of $(i)$).
Hence, we can assume that $I$ is well ordered and apply the method
of transfinite induction to construct a series of specialisations
$(P(K_{n,i}),P(L),\pi_{n,i})$ as
follows;\\

For $i=0$, set $(P(K_{n,0}),P(L),\pi_{n,0})=(P(K_{n}),P(L),\pi_{n})$\\

Given $i\in I$ with $i$ not a limit ordinal, let $v_{i+1}$ be the
next valuation in the enumeration. Let
$(K_{n}\{t\},\overline{v_{i+1}})$ be the completion of
$(K_{n}(t),v_{i+1})$ and let $\overline{v_{i+1}}$ also denote the
unique extension of this valuation to the algebraic closure
$K_{n}\{t\}^{alg}$. This defines a Krull valuation and hence a
specialisation $(P(K_{n}\{t\}^{alg}),P(L'),\pi_{n,i+1})$ where
$L'$ is the algebraic closure of the residue field of $v_{i+1}$,
having transcendence degree at most $1$ over $L$. Using arguments
as above, we can construct a specialisation $(P(L'),P(L),\pi)$.
Composing these specialisations, we obtain a specialisation
$(P(K_{n}\{t\}^{alg}),P(L),\pi_{n,i+1})$. (One can omit this step
by enumerating in $I$ only those valuations which preserve the
residue field $L$) Now, using Lemma 2.1 and Theorem 2.1,
amalgamate the specialisations
$(P(K_{n}\{t\}^{alg}),P(L),\pi_{n,i+1})$ and
$(P(K_{n,i}),P(L),\pi_{n,i})$ to form a specialisation;\\

$(P(K_{n,i}),P(L),\pi_{n,i})\prec (P(K_{n,i+1}),P(L),\pi_{n,i+1})$.\\

For $i$ a limit ordinal, we set;\\

$(P(K_{n,i}),P(L),\pi_{n,i})=\bigcup_{j<i}(P(K_{n,j}),P(L),\pi_{n,j})$\\

By the usual union of chains arguments we have that;\\

$(P(K_{n,j}),P(L),\pi_{n,j})\prec (P(K_{n,i}),P(L),\pi_{n,i})$ for
$j<i$.\\

Repeating this process, we obtain a structure $(P(K_{n+1}),P(L),\pi_{n+1})$ such that;\\

$(P(K_{n}),P(L),\pi_{n})\prec (P(K_{n+1}),P(L),\pi_{n+1})$.\\

It remains to check that this structure has the universal property
$(*)$ for $m=n+1$. Let $L_{n+1}$ be an algebraically closed
extension of $L$ with transcendence degree $n+1$ and
specialisation $(P(L_{n+1}),P(L),\pi)$. Let $v_{\pi}$ be the
corresponding valuation and its restriction to $L\subset
L_{n}\subset L_{n+1}$, a subfield of transcendence degree $n$. The
corresponding specialisation $(P(L_{n}),P(L),\pi)$ already factors
through $(P(K_{n}),P(L),\pi_{n})$ (\dag) and the valuation
$v_{\pi}$ appears as $v_{i}$ in the enumeration $I$ when
restricted to $L_{n}(t)$. By a standard result in valuation
theory, see \cite{Neu}, there exists an $L_{n}(t)$-embedding
$\tau:L_{n}(t)^{alg}\rightarrow L_{n}\{t\}^{alg}$ such that
$v_{\pi}=\overline{v_{i}}\circ\tau$ $(\dag\dag)$\ (see notation
above). Combining $(\dag)$ and $(\dag\dag)$, we obtain an
embedding $\alpha:(P(L_{n+1}),P(L))\rightarrow (P(K_{n,i}),P(L))$
such that $\pi=\pi_{n,i}\circ\alpha$. This
proves the result. It is now clear that the structure\\

$(P(K_{univ}),P(L),\pi_{univ})=\bigcup_{i>0}(P(K_{i}),P(L),\pi_{i})$\\

has the required universal property, is a model of
$T_{axioms}$ and;\\

$(P(K_{i}),P(L),\pi_{i})\prec (P(K_{univ}),P(L),\pi_{univ})$ for
$i>0$.\\
\end{proof}

\end{section}

\end{document}